\documentclass[12pt,letter]{article}

\usepackage{mystyle}
\usepackage{mythmstyle_report}

\begin{document}

\title{On the automorphism group of a Johnson graph}

\author{Ashwin~Ganesan%
  \thanks{Department of Electronics and Telecommunication Engineering, Vidyalankar Institute of Technology, Wadala, Mumbai, India. Email:  
\texttt{ashwin.ganesan@gmail.com}.}
}

\date{}  

\maketitle

\begin{abstract}
\noindent  The Johnson graph $J(n,i)$ is defined to the graph whose vertex set is the set of all $i$-element subsets of $\{1,\ldots,n\}$, and two vertices are joined whenever the cardinality of their intersection is equal to $i-1$.  In Ramras and Donovan [\emph{SIAM J. Discrete Math}, 25(1): 267-270, 2011], it is conjectured that if $n=2i$, then the automorphism group of the Johnson graph $J(n,i)$ is $S_n \times \langle T \rangle$, where $T$ is the complementation map $A \mapsto \{1,\ldots,n\} \setminus A$.  We resolve this conjecture in the affirmative.  The proof uses only elementary group theory and is based on an analysis of the clique structure of the graph.
\end{abstract}




\section{Introduction}

The Johnson graph $J(n,i)$ is defined to be the graph whose vertex set is the set of all $i$-element subsets of $\{1,\ldots,n\}$, and two vertices $A, B$ are said to be adjacent in this graph whenever $|A \cap B|=i-1$.  This graph has been well-studied in the literature (cf. \cite{Cuaresma:Giudici:Praeger:2008} \cite{Dabrowski:Moss:2000} \cite{Daven:Rodger:1999} \cite{Devillers:Giudici:Li:Praeger:2008} \cite{Etzion:Bitan:1996} \cite{Kim:Park:Sano:2010} \cite{Numata:1990} \cite{Ramras:Donovan:2011}).  The automorphism group of a graph is the set of all permutations of the vertex set of the graph that preserves adjacency \cite{Godsil:Royle:2001}.  In \cite{Ramras:Donovan:2011}, it is proved that if $n \ne 2i$, then the automorphism group of the Johnson graph $J(n,i)$ is isomorphic to $S_n$.  In \cite[Conjecture 1, p. 269]{Ramras:Donovan:2011} it is conjectured that if $n=2i$, then the automorphism group of $J(n,i)$ is isomorphic to $S_n \times \langle T \rangle$, where $T$ is the complementation map $A \mapsto 
A^c$ and $A^c:=\{1,\ldots,n\} \setminus A$.  In the present paper, this conjecture is resolved in the affirmative.  

Actually, the automorphism group of $J(n,i)$ for both the $n \ne 2i$ and $n=2i$ cases was already determined in \cite{Jones:2005}, but the proof given there uses heavy group-theoretic machinery.  The main result of   \cite{Ramras:Donovan:2011} was to provide a proof for the $n \ne 2i$ case that uses only elementary group theory; the proof is based on an analysis of the clique structure of the graph.  In \cite{Ramras:Donovan:2011} the authors leave the $n=2i$ case open but make a conjecture for this case.  We resolve this conjecture in the affirmative by providing a proof that again uses only elementary group theory and a similar analysis of the clique structure of the graph.

We first recall some basic facts about the Johnson graphs $J(n,i)$. Two vertices $A, B$ are adjacent in this graph iff their intersection $A \cap B$ has cardinality $i-1$, and this occurs exactly when the cardinality of their symmetric difference is 2.   The Johnson graph $J(n,i)$ is isomorphic to the Johnson graph $J(n,n-i)$; an explicit bijection between their vertex sets that preserves adjacency is the complementation map $T: A \mapsto A^c$.  Hence, without loss of generality we shall restrict our study of the Johnson graphs $J(n,i)$ to the case where $i \le n/2$.   Also, the graphs $J(n,1)$ are the complete graphs and hence are not very interesting.  The graphs $J(n,2)$ are the line graphs of complete graphs, and their automorphism groups are known.  Thus, when studying $J(n,i)$ henceforth, it is assumed that $i \ge 3$.  

Each permutation in $S_n$ acts in a natural way on the set of $i$-element subsets of $\{1,\ldots,n\}$, and this induced action on the vertices of $J(n,i)$ is an automorphism of the graph.  Also, distinct permutations in $S_n$ induce distinct permutations of the $i$-element subsets.  Hence $S_n$ is isomorphic to a subgroup of the automorphism group of $J(n,i)$.  In some cases, $S_n$ happens to be the (full) automorphism group of $J(n,i)$.  A special case of the results in \cite[Theorem 2(a)(c)]{Jones:2005} is that when $n \ne 2i$, the automorphism group of $J(n,i)$ is isomorphic to $S_n$; a special case of the results in \cite[Theorem 2(e)]{Jones:2005} is that when $n=2i$, the automorphism group of $J(n,i)$ is isomorphic to $S_n \times S_2$.  The proofs given in \cite{Jones:2005} use heavy group-theoretic machinery. An elementary combinatorial proof of the former result is given in \cite{Ramras:Donovan:2011}, and an elementary combinatorial proof of the latter result is given in the present paper.

The following is the main result proved in the present paper:

\begin{Theorem}
If $n=2i$, then the automorphism group of the Johnson graph $J(n,i)$ is $S_n \times \langle T \rangle$, where $T$ is the complementation map $A \mapsto A^c$. 
\end{Theorem}

For $\theta \in S_n$, let $\rho_\theta$ denote the permutation of the vertex set of $J(n,i)$ induced by $\theta$.  It is clear that $\{\rho_\theta: \theta \in S_n\}$ is a subgroup of the automorphism group of $J(n,i)$.  When $n=2i$, the subgroup $\langle T \rangle$ also acts as a group of automorphisms of $J(n,i)$:

\begin{Lemma}
 Suppose $n=2i$. Then the complementation map $T: A \mapsto A^c$ is an automorphism of the Johnson graph $J(n,i)$. 
\end{Lemma}

\noindent \emph{Proof}:
Let $A$ and $B$ be two vertices in $J(n,i)$.  We show that $A$ and $B$ are adjacent in $J(n,i)$ iff $A^c$ and $B^c$ are adjacent in $J(n,i)$. Recall that two vertices are adjacent in $J(n,i)$ iff their intersection has cardinality $i-1$.  The cardinality $|A^c \cap B^c| = n-|A \cup B| = n - (|A|+|B|-|A \cap B|) = n - 2i+|A \cap B|$, which equals $|A \cap B|$ since $n=2i$.  Since $A \cap B$ and $A^c \cap B^c$ have the same cardinality, the complementation map preserves adjacency and nonadjacency in $J(n,i)$.
\qed

The group $\{\rho_\theta: \theta \in S_n\} \langle T \rangle$ of automorphisms of $J(n,i)$ obtained so far can be expressed as a direct product:

\begin{Lemma}
 Let $T$ denote the complementation map $A \mapsto A^c$. The group $H:=\{\rho_\theta: \theta \in S_n\} \langle T \rangle$ of automorphisms of $J(2i,i)$ is isomorphic to the direct product $S_n \times \langle T \rangle \cong S_n \times S_2$.
\end{Lemma}

\noindent \emph{Proof}:  Observe that if $A$ is any $i$-element subset of $\{1,\ldots,n\}$, then $[\theta(A)]^c = \theta(A^c)$, whence $T$ and $\rho_\theta$ commute. It follows that $\{\rho_\theta: \theta \in S_n \} \langle T \rangle$ is a group and its two factors are normal subgroups. It remains to show that the two factors $\{\rho_\theta: \theta \in S_n\}$ and $\langle T \rangle$ have a trivial intersection.  By way of contradiction, suppose $T=\rho_\theta$ for some $\theta \in S_n$.  Then $\theta$ takes $\{1,\ldots,i-1,i\}$ to its complement $\{i+1,\ldots,2i\}$, and $\{1,\ldots,i-1,i+1\}$ to its complement $\{i,i+2,\ldots,2i\}$.  Hence $\theta$ takes the common elements $\{1,\ldots,i-1\}$ to $\{i+2,\ldots,2i\}$, and hence the remaining elements $\{i,i+1\}$ to $\{i,i+1\}$.  Take $A = \{2,\ldots,i-1,i,i+1\}$. Then $A^{\rho_\theta} \supseteq \{i,i+1\}$.  Thus $A^{\rho_\theta} \ne A^c$, a contradiction.  
\qed

\bigskip \noindent \textbf{Notation}. Fix a vertex $X$  of the graph $J(n,i)$.   Let $\mathcal{L}_i$ denote the set of vertices of $J(n,i)$ whose distance to $X$ is exactly $i$.  Thus, $\mathcal{L}_0 = \{X\}$, and $\mathcal{L}_1$ is the set $N(X)$ of neighbors of $X$.  Let $G$ denote the automorphism group of $J(n,i)$. The stabilizer of $X$ in $G$ is denoted $G_X$. 

We use the following additional notation from \cite{Ramras:Donovan:2011}. Each neighbor of a vertex $X$ in $J(n,i)$ is of the form $(X - \{p\}) \cup \{q\}$ for some $p \in X, q \notin X$. We denote this neighbor by $Y_{p,q}$.  For each $p \in X$, the set of neighbors $\{Y_{p,q}: q \notin X\}$ forms a clique, denoted by $\mathcal{Y}_p$.  The set $\{ \mathcal{Y}_p: p \in X \}$ is a partition of the set $N(X)$ of neighbors of $X$ into $i$ cliques, each of cardinality $n-i$.   Similarly, for each $q \notin X$, the set $\{Y_{p,q}: p \in X\}$ forms a clique, denoted by $\mathcal{Z}_q$. The set $\{\mathcal{Z}_q: q \notin X\}$ is a partition of $N(X)$ into $n-i$ cliques, each of cardinality $i$. Each maximal clique in $J(n,i)$ that contains the vertex $X$ is of the form  $ \{X\} \cup \mathcal{Y}_p$ for some $p \in X$ or of the form $\{X\} \cup \mathcal{Z}_q$ for some $q \notin X$ (cf. \cite[Lemma 1]{Ramras:Donovan:2011}).

We call each clique $\mathcal{Y}_p$ a clique of the first kind. Similarly, each clique $\mathcal{Z}_q$ is a clique of the second kind.  When $n \ne 2i$, the cardinality of a clique of the first kind is not equal to the cardinality of a clique of the second kind; thus, any automorphism of the graph that fixes the vertex $X$ must permute the set of cliques of the first kind in $N(X)$ amongst themselves.  On the other hand, when $n=2i$, the cliques in $N(X)$ of the first and second kind have the same cardinality, and so it is possible that there is an automorphism in $G_X$ that takes a clique of the first kind to a clique of the second kind. Indeed, we show below that such an automorphism exists and can be expressed in terms of the complementation map.

\begin{Proposition} \label{prop:G:acts:identically:as:rhothetaT} 
Suppose $n=2i$, and let $X$ be a vertex of $J(n,i)$ and let $g \in G_X$. Then there exist $\theta \in S_n$ and $i \in \{0,1\}$ such that the actions of $g$ and $\rho_\theta T^i$ on $\mathcal{L}_0 \cup \mathcal{L}_1$ are identical.
\end{Proposition}

\noindent \emph{Proof}:  Let $g \in G_X$. Then $g$ acts on the set $N(X)$ of neighbors of $X$, and hence permutes the maximal cliques in $N(X)$ amongst themselves.  Recall that these maximal cliques are either of the first kind or the second kind.  We consider two cases.  

First suppose that $g$ permutes the set of cliques in $N(X)$ of the first kind amongst themselves.  Since $g \in G_X$, $g$ acts bijectively on the set of all maximal cliques in $N(X)$, and so $g$ also permutes the set of cliques of the second kind amongst themselves. Hence $g: \mathcal{Y}_p \mapsto \mathcal{Y}_{\theta_1(p)}, \mathcal{Z}_q \mapsto \mathcal{Z}_{\theta_2(q)}$ for some $\theta_1 \in \Sym(X), \theta_2 \in \Sym(X^c)$.  Define $\theta \in S_n$ to be the map that takes $j$ to $\theta_1(j)$ if $j \in X$ and that takes $j$ to $\theta_2(j)$ if $j \in X^c$. As shown in \cite[p. 268]{Ramras:Donovan:2011}, the actions of $g$ and $\rho_\theta$ on $\mathcal{L}_0 \cup \mathcal{L}_1$ are identical.

For the rest of the proof, suppose that $g$ takes some clique of the first kind to a clique of the second kind.  So there exist $p' \in X, q' \notin X$ such that $g: \mathcal{Y}_{p'} \mapsto \mathcal{Z}_{q'}$.  We show that $g$ takes every clique of the first kind to some clique of the second kind.  Observe that $\mathcal{Z}_{q'}$ contains exactly one vertex from $\mathcal{Y}_p$, for each $p \in X$.  Any two cliques of the first kind are disjoint, and $g$ must map disjoint cliques to disjoint cliques.  Also, any two cliques of the second kind are disjoint, whereas a clique of the first kind and a clique of the second kind meet: $\mathcal{Y}_p \cap \mathcal{Z}_q \ne \phi$ since it contains $Y_{p,q}$. Thus, if $g$ takes a clique of the first kind to a clique of the second kind, then $g$ takes each clique of the first kind to some clique of the second kind.  Hence $g$ interchanges the set of cliques of the first kind and the set of cliques of the second kind.  

Thus $g: \mathcal{Y}_p \mapsto \mathcal{Z}_{\theta_1(p)}, \mathcal{Z}_q \mapsto \mathcal{Y}_{\theta_2(q)}$ for some $\theta_1: X \mapsto X^c$ and $\theta_2: X^c \mapsto X$.  Define $\theta \in S_n$ to be the map that takes $j$ to $\theta_1(j)$ if $j \in X$ and that takes $j$ to $\theta_2(j)$ if $j \in X^c$.  Recall that $\rho_\theta$ is defined as the action of $\theta$ induced on the vertex set of $J(n,i)$ and that $T$ denotes the complementation map $A \mapsto A^c$.  

We show that the actions of $g$ and $\rho_\theta T$ on $\mathcal{L}_0 \cup \mathcal{L}_1$ are identical.  It is clear that both the actions fix $\mathcal{L}_0 = \{X\}$.  For $g \in G_X$ implies $g$ fixes $X$.  And $X^{\rho_\theta T} = (X^c)^T = X$. Let $Y_{p,q}$ be a vertex in $\mathcal{L}_1$, and consider the action of $g$ and $\rho_\theta T$ on this vertex.  Recall that $Y_{p,q}$ is the unique vertex in the intersection $\mathcal{Y}_p \cap \mathcal{Z}_q$.  We have that $(\mathcal{Y}_p \cap \mathcal{Z}_q)^g = \mathcal{Z}_{\theta_1(p)} \cap \mathcal{Y}_{\theta_2(q)} = Y_{\theta_2(q), \theta_1(p)}$.   The vertex $Y_{p,q}$ has the same image under $\rho_\theta T$ as under $g$:
$(Y_{p,q})^{\rho_\theta T} = ((X-\{p\}) \cup \{q\})^{\rho_\theta T} = [(X^{\theta_1} - \{\theta_1(p)\}) \cup \{\theta_2(q)\}]^T = [(X^c - \{\theta_1(p)\}) \cup \{\theta_2(q)\}]^T = (X-\{\theta_2(q)\}) \cup \{\theta_1(p)\} = Y_{\theta_2(q),\theta_1(p)}$.  Thus, $g$ and $\rho_\theta T$ act identically on $\mathcal{L}_1$.
\qed

The following result, which is proved in \cite[Lemma 2 and Proposition 1]{Ramras:Donovan:2011}, does not use the condition that $n \ne 2i$ and hence also applies when $n = 2i$:

\begin{Lemma} \label{lemma:Le:equals:1}
 In the Johnson graph $J(n,i)$, if an automorphism $g$ fixes a vertex $X$ and each of its neighbors, then it is the trivial automorphism.
\end{Lemma}

We now complete the proof of the main theorem.

\begin{Corollary}
 If $n=2i$, then the automorphism group of the Johnson graph $J(n,i)$ is $S_n \times \langle T \rangle$, where $T$ is the complementation map $X \mapsto X^c$.
\end{Corollary}

\noindent \emph{Proof}:
Let $g \in G_X$.  By Proposition~\ref{prop:G:acts:identically:as:rhothetaT},  there exist $\theta \in S_n$ and $i \in \{0,1\}$ such that the action of $g$ and $\rho_\theta T^i$ are identical on $\mathcal{L}_0 \cup \mathcal{L}_1$.  Hence, $g^{-1} \rho_\theta T^i$ acts trivially on $\mathcal{L}_0 \cup \mathcal{L}_1$.  By Lemma~\ref{lemma:Le:equals:1}, $g^{-1} \rho_\theta T^i$ is the trivial automorphism of $J(n,i)$.  Hence $g = \rho_\theta T^i$.  This proves that every element in $G_X$ is one of the $2(i!)^2$ automorphisms specified in the proof above, i.e. every element in $G_X$ is either one of the $i! i!$ elements in $G_X$ that permutes the $i$ cliques of the first kind amongst themselves and the $i$ cliques of the second kind amongst themselves, or is one of the $i! i!$ elements in $G_X$ that interchanges the set of cliques of the first kind and the set of cliques of the second kind.  Hence $|G_X| = 2(i!)^2$. Finally, since the graph $J(n,i)$ is vertex-
transitive, the automorphism group $G$ has order $|G_X| {n \choose i} = 2(i!)^2 {n \choose i} = 2n!$.    Hence the group of automorphisms $S_n \times \langle T \rangle$ obtained above is the (full) automorphism group of $J(n,i)$.
\qed

 {
\bibliographystyle{plain}
\bibliography{refsaut}
}
\end{document}